\documentstyle [12pt] {article}

\title {Subalgebras with Converging Star Products
        in Deformation Quantization:\\
        An Algebraic Construction for $\complex \mbox{\LARGE P}^n$
        \vspace{1cm}}

\author {{\bf M. Bordemann\thanks{mbor@phyq1.physik.uni-freiburg.de}~,
            \addtocounter{footnote}{2}
            M. Brischle\thanks{brischle@phyq1.physik.uni-freiburg.de}~,
            \addtocounter{footnote}{-2}
            C. Emmrich\thanks{cemm@phyq1.physik.uni-freiburg.de}~,
             \addtocounter{footnote}{1}
            S. Waldmann\thanks{waldman@phyq1.physik.uni-freiburg.de}}\\[3mm]
             Fakult\"at f\"ur Physik\\Universit\"at Freiburg\\
           Hermann-Herder-Str. 3\\79104 Freiburg i.~Br., F.~R.~G\\[3mm]
       }

\date{FR-THEP-95/28 \\[1mm] December 1995}

\input{amssym.def}
\input{amssym}

\newcommand{\beq}{\begin{equation}}
\newcommand{\eeq}{\end{equation}}
\newcommand{\beqa}{\begin{eqnarray*}}
\newcommand{\eeqa}{\end{eqnarray*}}
\newcommand{\bea}{\begin{eqnarray}}
\newcommand{\eea}{\end{eqnarray}}
\newtheorem{lemm}{\sc Lemma }
\newtheorem{theo}{\sc Theorem }

\newtheorem{cor}{\sc Corollary }

\newcommand{\complex}{\mbox{C \hspace{-1.16em} \raisebox{-0.018em}{\sf l}}\;}
\newcommand{\qed}{\hfill $\Box $}
\newcommand{\pbl}{\pi^\ast}
\newcommand{\abla}{\partial z^{i_1}\cdots\partial z^{i_r}}
\newcommand{\ablb}{\partial \bar{z}^{i_1}\cdots\partial \bar{z}^{i_r}}

\newcommand{\rar}{\rightarrow}
\newcommand{\f}{\frac}
\newcommand{\p}{\partial}
\newcommand{\lb}{\label}
\newcommand {\CP} {\Bbb{C}\mbox{P}}
\newcommand {\CNull}{\Bbb{C}^{n+1}\setminus \{0\}}

\newcommand{\mustar} {*^\mu}
\newcommand {\newstar} {\, \tilde{*} \,}
\newcommand {\vp}{\varphi}

\begin{document}
\maketitle

\begin{abstract}
Based on a closed formula for a star product of Wick type on $\CP^n$,
which has been discovered in an earlier article of the
authors, we explicitly construct a subalgebra of the formal star-algebra
(with coefficients contained in the uniformly dense subspace of
representative functions with respect to the canonical action of the
unitary group) that consists of {\em converging} power series in the
formal parameter, thereby giving an elementary algebraic proof of a
convergence result already obtained by Cahen, Gutt, and Rawnsley. In
this subalgebra the formal parameter can be substituted by a real
number $\alpha$: the resulting associative algebras are
infinite-dimensional except for the case $\alpha=1/K$, $K$ a positive
integer, where they turn out to be isomorphic to the finite-dimensional
algebra of linear operators in the $K$th energy eigenspace of an
isotropic harmonic oscillator with $n+1$ degrees of freedom. Other
examples like the $2n$-torus and the Poincar\'e disk are discussed.
\end{abstract}

\vfill
\newpage

\section{Introduction}

The concept of deformation quantization as defined by
F.~Bayen, M.~Flato, C.~Fronsdal, A.~Lichnerowicz, and D.~Sternheimer
in 1978 (cf. \cite{bayen}) consists in a formal local deformation of the
commutative algebra (a so-called star product)
of all smooth complex-valued functions on a symplectic
manifold such that the first order commutator equals the Poisson
bracket and pointwise complex conjugation remains an antilinear involution.
The existence of star products on every symplectic manifold was proved by
M.~DeWilde and P.~B.~A.~Lecomte in 1983 (cf. \cite{DL 83}) and independently
by B.~Fedosov in 1985 \cite{Fed 85}, \cite{Fed 94}.
A third existence proof was given
by H.~Omori, Y.~Maeda and A.~Yoshioka \cite {OMY 91}.

One of the problems with these star products is the fact that the formal
series involved is shown to never converge on the space of {\em all}
complex-valued smooth functions (see e.g.\cite{Rub 84}),
i.e. for every complex
number there are two functions whose star product diverges when the
formal parameter is substituted for that number.

This paper is a continuation of our work \cite{uns} in which we gave
a closed formula for a star product of Wick type on complex projective
space by a version of quantum phase space reduction.
Star products on complex projective space have already been
constructed by
H.~Omori, Y.~Maeda and A.~Yoshioka \cite {OMY 93} and by
C.~Moreno \cite {Mor 86} in a less explicit way.
$\CP^n$ can also be regarded as a coadjoint orbit of the unitary
group $U(n+1)$ (see the work of D.~Arnal, J.~Ludwig and M.~Masmoudi
for the existence of covariant star products on more general coadjoint
orbits in \cite {ALM 94} and references therein).

The aim of the present paper is twofold: firstly, we use our formula to
explicitly compute the star product
for a uniformly dense subalgebra of representative functions for
the $U(n+1)$-action on complex projective space. Secondly, we would
like to use the example $\CP^n$ in order to illustrate
the following {\em algebraic} procedure to deal with the above convergence
problem:
\begin{enumerate}
\item Find a complex subalgebra ${\cal U}$ of the $*$-algebra of
      all formal power series in $\nu$ with smooth coefficients
      such that  1.) all elements of ${\cal U}$ are power series with
      infinite radius of convergence in $\nu$ and 2.) their coefficients
      may be choosen in a ``sufficiently large'' (e.g. uniformly
      dense) subspace of the space of all smooth complex-valued functions.
\item Verify whether the subspace ${\cal I}_{\alpha}$ of ${\cal U}$
      defined by
      \beq
          {\cal I}_{\alpha} := \{\Phi(\nu) \in {\cal U} | \Phi(\alpha)=0 \}
      \eeq
      is a star-ideal of ${\cal U}$.
\item Identify the quotient ${\cal U}/{\cal I}_{\alpha}$ with the
      associative algebra of quantum observables related to
      the ``$\hbar$-value'' $\alpha$ and try to find a
      representation of this quotient in some Hilbert space.
\end{enumerate}
{}From the physical point of view it is often required that the subalgebra
$\cal U$ contains certain ``important observables'' which are in some cases
related to the presence of additional symmetries of the classical phase
space. For a general symplectic manifold the viability of the above
procedure (in particular the existence of a sufficiently large subalgebra
$\cal U$) is --to our best knowledge-- an open problem in the theory of
deformation quantization.

Nevertheless, in several examples having a large symmetry group
(e.g. the Moyal product on the 2-torus or the Wick product on
$\Bbb{C}^{n+1}$) the above programme may be carried out.
For such manifolds the space of representative functions
of the symmetry group plays a prominent r\^ole for the construction
of $\cal U$: in a remarkable article of M.~Cahen, S.~Gutt, and J.~Rawnsley
\cite{CGR II} the convergence of a star product for these functions has
been proved for all compact Hermitean spaces (in particular for $\CP^n$)
by analytic methods of complex differential geometry.
They start from the
finite-dimensional operator algebras of geometric quantization
in tensor powers of a very ample regular prequantum line bundle over
a compact K\"ahler manifold and use coherent states (see \cite{Ber 74},
\cite{Raw 77}) to first construct star products for the Berezin-Rawnsley
symbols (\cite{Ber 74}, \cite{Raw 77}) for each tensor power separately.
In a second step
an asymptotic expansion of these star products in the inverse tensor
power is shown to define a local star product on the manifold where
the formal parameter appears as a sort of interpolation of the inverse
tensor powers.

The approach of this paper is in some sense reverse
to the programme of M.~Cahen, S.~Gutt, and J.~Rawnsley
(cf. \cite{CGR I}, \cite{CGR II},\cite{CGR III}, \cite{CGR IV})
and only makes use of elementary algebraic methods:
We start from the explicit star product on $\CP^n$ (see \cite{uns})
and define $\cal U$ as a certain proper subspace of the space of all
polynomials in the formal parameter with coefficients in the
uniformly dense subspace of representative functions
for the unitary group $U(n+1)$. Since all occurring star products
can explicitly be computed the analysis of the ideals ${\cal I}_\alpha$
and the quotient algebras ${\cal U}/{\cal I}_\alpha$
becomes relatively simple. The main result is that for inverse integer
values of
the formal parameter $\nu$ the quotient algebras turn out to be
finite-dimensional full complex matrix-algebras whereas all
noninteger quotients
are of infinite dimension and define converging star products on the
space of representative functions for these values as has been stated in
\cite{CGR II}.

{}From the physical point of view one can
regard the finite-dimensional quantum algebras as the set of all
quantum observables restricted to the eigenspace of integral energy of a
harmonic oscillator of $n+1$ degrees of freedom (where the ground state
energy is zero for the Wick quantization rule).

The advantage of the above algebraic programme of prescribing a real value
to the formal parameter of deformation quantization is that one may hope
to transfer it to physical situations with an infinite number of degrees
of freedom, i.e. field theories, where the powerful analytical methods
in the theory of finite-dimensional manifolds are no longer well-defined.

The paper is organised as follows: after briefly reviewing the concepts
and formulas of our last paper \cite{uns} in section 2 we compute the
star product of two Berezin-Rawnsley symbols in section 3 and discuss
the unitary symmetry action on the star product.
Section 4 is then devoted to define the algebra ${\cal U}$, to compute the
ideals ${\cal I}_{\alpha}$ and the quotient algebras
${\cal U}/{\cal I}_{\alpha}$.
In section 5 we briefly consider other phase spaces already dealt with
in the literature for which the above
programme works: the complex vector space $\Bbb{C}^{n+1}$ with the
Wick product, the $2n$-torus $T^{2n}$ with the Moyal product, and the
Poincar\'e disk. For this last example we can recover the formula of the
star product for the corresponding Berezin-Rawnsley symbols given by
Cahen, Gutt, and Rawnsley in \cite{CGR III}.

\vspace {0.3cm}
\noindent

{\it Notation:}
Throughout this paper we use the Einstein summation convention,
i.e. summation over repeated indices is automatic. Moreover, the symbol
$F(z)$ for a complex-valued function $F$ of a complex vector $z$ does
{\em not} necessarily imply that $F$ is holomorphic.

\section{Review of star products on complex projective space}

In this section we shall give a short review of earlier work \cite{uns}
in which we derived an explicit formula for a star product
of Wick-type on the complex projective space $\CP^n$.

Let $\pi: \CNull \to \CP^n$ be the canonical projection of a complex
vector $z$ onto the complex ray through it. Let $x:=\bar{z}^iz^i$.
The usual Wick product on $\CNull$ of two complex-valued functions
$F, G \in C^\infty(\CNull)$ is given
as the following formal power series in the parameter $\lambda$:
\beq \lb{Wickusual}
    F\ast G = \sum_{r=0}^\infty
    \frac{\lambda^r}{r!}\frac{\partial^r F}{\abla}
    \frac{\partial^rG}{\ablb}.
\eeq
We have called a function $F \in C^\infty (\CNull)$ {\em homogeneous}
iff it is invariant under the natural action of the group
$\Bbb{C}\setminus \{0\}$. These functions are precisely given by
pull backs of functions $f \in C^\infty(\CP^n)$,
$F = f \circ \pi = \pi^*f$. We have called
a function $R \in C^\infty(\CNull)$ {\em radial} iff it is a function
of $x$, i.e. iff there is a smooth function $\rho: \Bbb{R}^+ \to \Bbb{C}$
such that $R = \rho \circ x$.
We have defined a formal differential operator
$S: C^\infty(\CNull)[[\lambda]] \to C^\infty(\CNull)[[\lambda]]$
depending only on $x$ and
$\p_x := \f{1}{2x}(z^i\p_{z^i}+\bar{z}^i\p_{\bar{z}^i})$
whose standard symbol
$\hat{S}(x,\alpha) := (Se_\alpha)(x) e^{-\alpha x}$ ($e_\alpha$ denoting
the exponential function $x \mapsto e^{\alpha x}$ for $\alpha\in\Bbb{R}$)
is given by (setting the series $D$ in \cite [eqn. 9] {uns} equal to $1$):
\beq \lb{Ssymb}
    \hat{S}(x,\alpha) = \exp\left(\frac{x}{\lambda}
    \log(1+\lambda\alpha) - \lambda\alpha)\right).
\eeq
$S$ and its inverse $S^{-1}$ act trivially on homogeneous functions, i.e.
$SF = F = S^{-1}F$, but do in general transform radial functions
into radial ones, in particular (for a nonnegative integer $r$)
\cite [eqn. 14] {uns}
\beq \lb{xrel}
    Sx=x,
    ~~~~
    Sx^r = x^r\prod_{k=0}^r \left(1 - k\frac{\lambda}{x}\right),
    ~~~~
    Sx^{-r} = x^{-r}\prod_{k=0}^r \left(1 + k\frac{\lambda}{x}\right)^{-1}.
\eeq
We have used $S$ as an equivalence transformation for a modified
Wick product of two functions $F,G\in C^\infty(\CNull)$:
\beq \lb{Wickmodif}
   F \newstar G := S (S^{-1}F \ast S^{-1}G).
\eeq
For two radial functions $R_1,R_2$
and a homogeneous function $F$ on $\CNull$ this new star product is
just pointwise multiplication \cite [eqn. 12] {uns}:
\beq \lb{radradhom}
    R_1 \newstar R_2 = R_1R_2 = R_2 \newstar R_1,
    ~~~~
    R_1 \newstar F = R_1F = F \newstar R_1,
\eeq
whereas for two smooth homogeneous functions $F, G$ we get
\beq \lb{homhom}
  (F \newstar G) (z)
  = \sum_{r=0}^\infty\frac{1}{r!}
  \left(\frac{\lambda}{x}\right)^r
  \prod_{k=1}^r \left( 1+k\frac{\lambda}{x} \right)^{-1}
  x^r\frac{\partial^r F}{\abla}(z)\frac{\partial^r G}{\ablb}(z).
\eeq
The main result of \cite{uns} was the fact that this formula can
directly be projected to $\CP^n$ by phase space reduction via
the $U(1)$-momentum map $J: \CNull \to \Bbb{R}: z \mapsto -\f{x}{2}$
of the canonical $U(1)$-action on $\CNull$: for a negative real
number $\mu$ and a $U(1)$ invariant function $F$ in
$C^\infty(\CNull)$ we write $F_\mu$ for the unique function in
$C^\infty(\CP^n)$ obtained by first restricting $F$ to the odd
sphere $J^{-1}(\mu)$ and then projecting it to $\CP^n$
\cite [4.3] {AM 85}. Then the formula
\beq \lb{reduct}
    F_\mu \mustar G_\mu:=(F \newstar G)_\mu
\eeq
was shown to define a star product on $\CP^n$. The explicit form
of $\mustar$ is obtained by replacing $\lambda/x$ by $\lambda/(-2\mu)$
in (\ref{homhom}) and noting that the bidifferential operator
$\tilde{M}_r(f,g)(\pi(z)) := x^r\f{\p^r \pi^*f(z)} {\abla}
\f{\p^r \pi^*g(z)} {\ablb}$
is well-defined on $f, g \in C^\infty(\CP^n)$. For simplicity we shall
work with the redefined formal parameter $\nu:=\lambda/(-2\mu)$ in what
follows.

\begin{lemm} :
   The standard symbol of $S^{-1}$ is described by the formula
    \beq \label {zeitent}
      \widehat{S^{-1}}(x,\alpha)
      = e^{\f{x}{\lambda}(e^{\alpha\lambda}-1-\alpha\lambda)}
      = e^{\ast\alpha x}e^{-\alpha x}
    \eeq
   where the last term involves the star-exponential \cite{bayen} of
   the function $x$ with respect to the usual Wick product
   (\ref{Wickusual}).
\end{lemm}
{\sc Proof}: Since $Se_\beta=e_{\f{1}{\lambda}\log(1+\beta\lambda)}$
we obviously get $e_\beta=S^{-1}e_{\f{1}{\lambda}\log(1+\beta\lambda)}$ which
proves the first equation after the substitution
$\alpha:=\f{1}{\lambda}\log(1+\beta\lambda)$. Secondly, note that
\[
  e^{\ast\alpha x} = S^{-1}Se^{\ast\alpha S^{-1}x} =
  S^{-1}e^{\newstar\alpha x}
  \stackrel{(\ref{radradhom})} {=} S^{-1}e^{\alpha x} =
  \widehat{S^{-1}}(x,\alpha)e^{\alpha x}
\]
which proves the second equality. \qed

Remark: Note that the function $H := \f{1}{2} x$ equals the usual
Hamiltonian function of an isotropic harmonic oscillator in $n+1$
degrees of freedom. The above star-exponential of $x$ for
$\alpha = - it / 2\hbar$ and $\lambda = 2\hbar$
then corresponds to the quantum mechanical time development
operator for this system.

\section{A star product for representative functions on
         complex projective space}

Let $p: \Bbb{C}^{n+1} \to \Bbb{C}: z \mapsto P(z)$ be a polynomial
function (in the $2n+2$ variables
$z^0, \ldots, z^n$, $\bar{z}^0, \ldots, \bar{z}^n$).
We shall call $p$ {\em homogeneous of degree $(k,k)$} for a nonnegative
integer $k$ iff $p(\lambda z) = (\lambda\bar{\lambda})^k p(z)$
for all $\lambda\in\Bbb{C}\setminus\{0\}$.
We denote by ${\cal E}_k$ the following subspace of $C^\infty (\Bbb{C} P^n)$:
\beq \label{dieEks}
   {\cal E}_k := \left\{\phi \left|
   \mbox{ there is a homogeneous polynomial $p_k$ of degree $(k,k)$ s.t.~}
           (\pi^\ast\phi)(z)=\frac{1}{x^k}p_k(z) \right\}\right..
\eeq

\begin{lemm} :
\begin{enumerate}
\item For each integer $k \geq 0$: ${\cal E}_k\subset{\cal E}_{k+1}$.
\item ${\cal E}:=\bigcup_{k=0}^\infty{\cal E}_k$ is a filtered
      subalgebra of $C^\infty(\Bbb{C} P^n)$
      with respect to the pointwise multiplication. It is closed under
      complex conjugation.
\item $\cal E$ separates points and is therefore a dense subalgebra
      of $C^\infty(\Bbb{C}P^n)$ with respect to the uniform topology.
\end{enumerate}
\end{lemm}
{\sc Proof}:
i)   Clearly $\f{p_k}{x^k}=\f{xp_k}{x^{k+1}} \in \pi^* {\cal E}_{k+1}$.
ii)  This is obvious.
iii) Consider the complex rays $\pi(z_{(1)})\not=\pi(z_{(2)})$ for
$z_{(1)},z_{(2)}\in\Bbb{C}^{n+1}\setminus 0$. This is true iff
$z_{(1)},z_{(2)}$ are linearly independent iff there is a
$y\in\Bbb{C}^{n+1}$ such that $\langle y,z_{(1)}\rangle = 1$
and $\langle y,z_{(2)}\rangle = 0$ where $\langle y,z \rangle$ denotes
the standard sesquilinear form $\langle y,z\rangle=\bar{y}^kz^k$. Then
$\phi\in {\cal E}_1$ defined by $\phi(\pi(z)):=\frac{|<y,z>|^2}{x}$
separates $\pi(z_{(1)})$ and $\pi(z_{(2)})$. The density of $\cal E$
follows from the Stone-Weierstrass Theorem. \qed \\

Consider the standard action of the unitary group $U(n+1)$ on $\CNull$:
$(g,z)\mapsto gz=:\Phi_g(z)$ and its induced action on $\CP^n$:
$(g,\pi(z))\mapsto \pi(gz)=:\Psi_g(\pi(z))$. A smooth
complex-valued function $f$ on $\CP^n$ is called {\em representative} with
respect to the $U(n+1)$ action iff
\beq \label{rep}
    \Bbb{C}{\rm -span} \{f \circ \Psi_g | g\in U(n+1)\}
    \quad \mbox{ is finite-dimensional.}
\eeq
We now get a characterization of $\cal E$ which should be fairly standard
(see \cite[sec. 3, Lemma 1]{CGR II}):
\begin{lemm} \lb{Eistrep}:
  The algebra $\cal E$ is equal to the set of all representative
  functions on $\CP^n$.
\end{lemm}
{\sc Proof}:
Since $x$ and the finite-dimensional space of all homogeneous
polynomials of degree $(k,k)$ are invariant under $U(n+1)$
it follows that $\cal E$ consists of representative functions.
In order to prove the reversed inclusion we can use a more general argument:
$\CP^n$ is a homogeneous space $G/H$ for the compact Lie group
$G=U(n+1)$ with compact isotropy subgroup $H=U(1)\times U(n)$. Now the
space of all representative functions $f$ on $G/H$ (defined as in (\ref{rep})
with $U(n+1)$ replaced by any compact Lie group $G$) is clearly in
one-to-one correspondence with the space of its pull-back
to $G$ under the natural projection $G\rar G/H$: this in turn is given by
the space of all $H$ right invariant representative functions on $G$ with
respect to left multiplication. Both this space and the pull-back of $\cal E$
to $G$ are $G$-modules and are therefore closed in the space of all
representative functions on $G$ with respect to the uniform topology on $G$
(see e.g. \cite{BtD 85}, p.126, Prop.(1.4) (iii)). Since the pull-back
obviously is a continuous closed linear map with respect to the uniform
topologies
and $\cal E$ is dense in the space of all representative functions on
$G/H$ thanks to the previous Lemma it follows that $\cal E$ is equal to that
space.
\qed \\

Yet another equivalent description of $\cal E$ is obtained in terms
of the {\em Berezin-Rawnsley symbols} known from geometric quantization
(\cite{Ber 74}, \cite{Raw 77}, \cite{CGR I}): for a fixed nonnegative
integer $k$ take the vector space ${\cal H}^{(k)}$ of complex-valued
holomorphic polynomials $\psi$ on $\Bbb{C}^{n+1}$, which are
homogeneous of degree $k$, i.e. $\psi(\lambda z)=\lambda^k\psi(z)$ for
all complex numbers $\lambda$.
The dimension of this space is clearly $N:={n+k \choose k}$ and we shall
henceforth identify ${\cal H}^{(k)}$ with $\Bbb{C}^N$ by means of the base
$(z^{i_1}\cdots z^{i_k})$. Consider
the space $B({\cal H}^{(k)})$ of complex linear endomorphisms of
${\cal H}^{(k)}$. Any $A\in B({\cal H}^{(k)})$ can be viewed as
a complex $N\times N$ matrix
$A_{i_1\cdots i_k,j_1\cdots j_k}$ where each of the indices
$i_1,\ldots,i_k,j_1,\ldots,j_k$ ranges over $0,1,\ldots,n$ and the
matrix elements are symmetric with respect to all permutations among
the $i_1,\ldots, i_k$ and among the $j_1,\ldots, j_k$. To each
$A\in B({\cal H}^{(k)})$ one can associate
the polynomial function
\beq \label{sigmaoben}
   \tilde{\sigma}(A):\Bbb{C}^{n+1}\rar\Bbb{C}:z\mapsto
      \bar{z}^{i_1}\cdots\bar{z}^{i_k}z^{j_1}\cdots z^{j_k}
                 A_{i_1\cdots i_k,j_1\cdots j_k}.
\eeq
Clearly, $\tilde{\sigma}(A)$ is homogeneous of degree $(k,k)$, and
by counting dimensions it can be seen that every homogeneous polynomial
function of degree $(k,k)$ is of that form.
The {\em Berezin-Rawnsley symbol $\sigma (A)$ associated to $A$}
is then defined by
\beq \label{BerRaw}
    \sigma(A) : \CP^n \to \Bbb{C}:
    \pi(z) \mapsto \f{\tilde{\sigma}(A)(z)}{x^k}.
\eeq
Then the following corollary is clear:
\begin{cor} :
    For each nonnegative integer $k$ the space ${\cal E}_k$ is spanned
    by all the Berezin-Rawnsley symbols $\sigma(A)$ associated to
    $A\in B({\cal H}^{(k)})$.
\end{cor}

Both on $\CNull$ and $\CP^n$ there are momentum maps (see \cite{AM 85} for
general definitions) for the $U(n+1)$-action which can be expressed in
terms of Berezin-Rawnsley symbols contained in ${\cal E}_1$:
\begin{lemm} : Let $\goth{u} (n+1)$ denote the space of
    complex antihermitean $(n+1)\times (n+1)$-matrices,
    i.e the Lie algebra of the unitary group $U(n+1)$. Then:
    \begin{enumerate}
    \item[i)] The following map
              \beq
                   \tilde{P}:\CNull\to \goth{u}(n+1)^* :
                   z \mapsto (A \mapsto \tilde{\sigma}(A))
              \eeq
              is a momentum map for the $U(n+1)$-action on $\CNull$.
    \item[ii)] The following map
               \beq
                   P:=\tilde{P}_\mu:\CP^n \to \goth{u}(n+1)^*:
                   \pi(z) \mapsto (A\mapsto -2\mu \sigma(A)(\pi(z))
                   = \tilde{\sigma}(A)_\mu (\pi(z))~)
               \eeq
               is a momentum map for the $U(n+1)$-action on $\CP^n$.
   \end{enumerate}
   Recall that the index $\mu$ refers to the $U(1)$-phase space reduction
   of $\CNull$ by the $U(1)$ momentum map $J(z)=-\f{x}{2}$.
\end{lemm}
{\sc Proof}:
i) The Hamiltonian vector field of $P(A)$ is given by
$X_{\tilde{\sigma}(A)}(z)$ $=$
$\frac{2}{i} (A_{ij}z^j\partial_{z^i}
 + \bar{A}_{ij}\bar{z}^j\partial_{\bar{z}^i})$,
since $A_{ij}=-\overline{A}_{ji}$ which obviously equals the infinitesimal
generator $A_{\CNull}$ of the $U(n+1)$-action.
The $Ad^\ast$-equivariance of this map is obvious.\\
ii) Recall the projection $\pi_\mu:J^{-1}(\mu)\rar \CP^n$ to the reduced
phase space (compare \cite[4.3]{AM 85}). Since
$\Phi_g\circ\pi_\mu=\pi_\mu\circ\Psi_g$ it is clear that
$T\pi_\mu A_{\CNull}=A_{\Bbb{C}P^n}\pi_\mu$ where $A_{\Bbb{C}P^n}$ is the
infinitesimal generator of the $U(n+1)$-action on $\CP^n$. Using the
identity
$A_{\Bbb{C}P^n}\pi_\mu(z)=X_{\tilde{\sigma}(A)_\mu}\pi(z)$ we get by
phase-space reduction
\[
  T_z\pi_\mu X_{\tilde{\sigma}(A)_\mu}(z)=
  T_z\pi_\mu A_{\CNull} (z)=
  X_{\tilde{\sigma}(A)_\mu}(\pi(z)).
\]
The $Ad^\ast$-equivariance of this map follows at once from
the $Ad^\ast$-equivariance of $\tilde{P}$. \qed\\

We shall now compute star products for elements of $\cal E$:
Set $\nu := \frac{\lambda}{-2\mu}$, $\nu^{(0)} := 1$, $\nu^{(1)} := 1$
and
\beq \label {nudef}
    \nu^{(k)} := (1-\nu) \dots (1-(k-1)\nu)
\eeq
\begin{theo}\label{one} :
For $f\in{\cal E}_k$, $g\in{\cal E}_l$
we get
\beq \lb{oneeins}
    (\pbl f)\newstar (\pbl g)(z)=\sum_{r=0}^{min(k,l)}
    \frac{\lambda^r}{r!}\frac{S(x^{k+l-r})}{(Sx^k)(Sx^l)}\frac{1}{x^{k+l-r}}
    \frac{\partial^r(x^k\pbl f)}{\abla}(z)
    \frac{\partial^r(x^l\pbl g)}{\ablb}(z)
\eeq
and
\beq \label {StarEqn}
    f \mustar g = \sum_{r=0}^{min(k,l)}
    \frac{\nu^r}{r!}\frac{\nu^{(k+l-r)}}{\nu^{(k)}\nu^{(l)}}
    M_r^{(k,l)}(f,g)
\eeq
with the following bidifferential operator on $C^\infty(\CP^n)$
\beq \lb{dieHs}
    \pi^\ast M^{(k,l)}_r(f,g)
    := \frac{1}{x^{k+l-r}}
    \frac{\partial^r(x^k\pi^\ast f)}{\abla}
    \frac{\partial^r(x^l\pi^\ast g)}{\ablb}
\eeq
\end{theo}
{\sc Proof} : The usual Wick product (\ref{Wickusual}) of the two
{\em polynomials}
$\tilde{f_k}:=x^k\pbl f,\tilde{g_l}:=x^l\pbl g$ is given by
\beq
     \tilde{f}_k\ast\tilde{g}_l =
     \sum_{r=0}^{min(k,l)} \frac{\lambda^r}{r!}
     \frac{\partial^r\tilde{f}_k} {\abla}
     \frac{\partial^r\tilde{g}_l}{\ablb}.
\eeq
Using definition (\ref{Wickmodif}) and the formulas
(\ref{radradhom}) and (\ref{xrel}) we find
\beqa {(Sx^k)(Sx^l)(\pbl f)\tilde{\ast} (\pbl g) }
 & = & \left( S(x^k\pbl f)\right)\tilde{\ast}
       \left( S(x^l\pbl g) \right) = S((x^k\pbl f)\ast (x^l\pbl g)) \\
 & = & S\left( \sum_{r=0}^{min(k,l)}\frac{\lambda^r}{r!}
       \frac{\partial^r\tilde{f}_k}{\abla}
       \frac{\partial^r\tilde{g}_l}{\ablb}
       \right) \\
 & = & \sum_{r=0}^{min(k,l)}\frac{\lambda^r}{r!}
       S\left( x^{k+l-r}\frac{1}{x^{k+l-r}}
       \frac{\partial^r\tilde{f}_k}{\abla}
       \frac{\partial^r\tilde{g}_l}{\ablb}
       \right) \\
 & = & \sum_{r=0}^{min(k,l)}\frac{\lambda^r}{r!}
       S\left(x^{k+l-r} \right) \frac{1}{x^{k+l-r}}
       \frac{\partial^r\tilde{f}_k}{\abla}
       \frac{\partial^r\tilde{g}_l}{\ablb}
\eeqa
which proves the first equation. The second equation immediately follows
by the reduction formula (\ref{reduct}).
\qed \\

Remark: Thanks to the straight forward recursion formulas
$M^{(k+1,l)}_r = M^{(k,l)}_r+r(l-(r-1)) M^{(k,l)}_{r-1}$ and
$M^{(k,l+1)}_r = M^{(k,l)}_r+r(k-(r-1)) M^{(k,l)}_{r-1}$ it can
easily be checked that the right hand side
of (\ref{oneeins}) is well-defined, i.e. if e.g. $f$ is regarded as an
element of ${\cal E}_{k+a}$ for a positive integer $a$.

\begin{cor}\lb{COR}: For $A, B \in B({\cal H}^{(1)})$ we have
\beq
    \sigma(A)^k \mustar \sigma(B)^l =
    \sum_{r=0}^{min(k,l)}
    \frac{\nu^r}{r!}\frac{k! \: l!} {(k-r)!(l-r)!}
    \frac{\nu^{(k+l-r)}}{\nu^{(k)}\nu^{(l)}}
    \sigma(AB)^r\sigma(A)^{k-r}\sigma(B)^{l-r}
\eeq
\end{cor}
{\sc Proof}: This is easily seen by setting $f = \sigma(A)^k$ and
$g = \sigma(B)^l$ and using Theorem \ref{one}.
\qed \\

We shall now show that the momentum map $P$ for the $U(n+1)$-action
on $\CP^n$ is even a quantum momentum map. More precisely:
\begin{lemm} \lb{quantmom}:
   $P$ is a quantum momentum mapping on $\CP^n$ for the
   $U(n+1)$-action, i.e.
   for every smooth function $\phi: \CP^n \to \Bbb C$ the
   following equation holds:
   \beq \label {invarianz}
        P(A) \mustar \phi - \phi \mustar P(A)
        = \frac{i\lambda}{2}\{P(A), \phi\}_\mu ,
   \eeq
   i.e. the star product $\mustar$ is strongly $U(n+1)$-invariant.
\end{lemm}
{\sc Proof}:
First we prove the equation
$\pi^*\sigma(A) \newstar f - f \newstar \pi^*\sigma(A) = \frac{i\lambda}{2}
\{\pi^*\sigma(A), f\}$ for $f = \pi^*\phi$. Then eqn. (\ref{invarianz})
follows by eqn.(\ref{reduct}).
We have the strong invariance of the Wick product
\[
   \tilde{\sigma}(A) * f - f * \tilde{\sigma} (A) =
   \frac{i\lambda}{2} \{\tilde\sigma (A), f \}.
\]
With the equivalence transformation $S$ we find
\[
   S^{-1} \Big( S\tilde\sigma (A) \newstar Sf -
                 Sf \newstar S\tilde\sigma (A) \Big)
   = \frac{i\lambda}{2} \{\tilde\sigma (A), f \}
\]
and with $\pi^*\sigma (A) = x \tilde \sigma (A)$ and $Sf = f$
this leads to
\[
   Sx \Big(
   \pi^*\sigma (A) \newstar f - f \newstar \pi^*\sigma (A) \Big)
   =
   \frac{i\lambda}{2} S \left(\{x \pi^*\sigma (A), f \}\right).
\]
The Poisson bracket $\{x\pi^*\sigma (A), f\} = x \{\pi^*\sigma (A), f\}$ is
again homogeneous. Hence the right hand side of the last equation is simply
$x \{\pi^*\sigma(A), f\}$. With $Sx = x$ the proof is complete.
\qed

Remark: The case $k=l=1$ in Cor.2 (\ref{COR}) shows that the functions
$\sigma(A)$ and $\sigma(B)$ commute with respect to $\mustar$ iff
the operators $A$ and $B$ commute. Let $A_1,\ldots, A_n$ be a linearly
independent set of commuting traceless Hermitean matrices. It follows in
particular that the functions $\sigma(A_1), \ldots, \sigma(A_n)$
are functionally
independent and in involution, i.e. the Poisson bracket of $\sigma(A_i)$
with $\sigma(A_j)$ vanishes for all $i,j$. In other words they form a
{\em completely integrable system} in the sense of Liouville on $\CP^n$.
Note that this system is nontrivial in the sense that there is no global
chart of action-angle-variables, i.e. $\CP^n$ is not symplectomorphic to
some $T^r \times \Bbb{R}^{2n-r}$ for a nonnegative integer $r\leq n$. The
above Corollary now
implies that these functions also commute with respect to the star
product and can be viewed as a {\em quantum integrable system} on
$\CP^n$.

\section{Construction of the subalgebra $\cal U$, the ideals
        $\cal I_\alpha$, and the quotients ${\cal U}/{\cal I}_\alpha$}

If we take two functions $\phi_k \in {\cal E}_k$ and
$\psi_l \in {\cal E}_l$ and multiply them with the $\nu$-polynomial
$\nu^{(k)}$ and $\nu^{(l)}$, respectively, then formula
(\ref{StarEqn}) shows that their $\mustar$-product
contains only {\em polynomials} in the parameter $\nu$.
Moreover, that
particular combination is restored, i.e. the functions
$M_t^{(k,l)} (\phi_k, \psi_l) \in {\cal E}_{k+l-t}$ appear only
in combination with $\nu^{(k+l-t)}$.
This motivates the definition of the following subspaces of
$C^\infty(\CP^n)[[\nu]]$:
\beq \lb{dieUks}
    {\cal U}_0 := {\cal E}_0 \mbox{ and }
    {\cal U}_k := \nu^{(k)}{\cal E}_k + \nu\nu^{(k-1)} {\cal E}_{k-1}
                  + \cdots + \nu^{k-1}\nu^{(1)}{\cal E}_1 +
                  \nu^k{\cal E}_0 \qquad \forall~k \in \Bbb{N}
\eeq
where each ${\cal U}_k$ is a (not necessarily direct) sum of subspaces of
$C^\infty(\CP^n)[[\nu]]$.
We denote by ${\cal E}[\nu]$ the {\em polynomials} in $\nu$ with
coefficients in ${\cal E}$.
The following theorem describes the structure of ${\cal U}$.
\begin{theo} :
\begin{enumerate}
\item For each integer $k$: ${\cal U}_k\subset {\cal U}_{k+1}$.
      Define ${\cal U}:=\bigcup_k^\infty {\cal U}_k$. Then
      $\cal U$ is a proper $\Bbb{C}[\nu]$-submodule of ${\cal E}[\nu]$.
\item $\cal U$ is a filtered subalgebra of
      $C^\infty(\Bbb{C}P^n)\left[ [\nu ] \right]$
      with respect to $\mustar$, i.e.
      \[ {\cal U}_k \mustar {\cal U}_l\subset {\cal U}_{k+l} \]
\end{enumerate}
\end{theo}
{\sc Proof :}
{\it i.)} Let $\Phi\in {\cal U}_k$. $\Phi$ is of the form
$\Phi=\sum_{r=0}^k\nu^r\nu^{(k-r)}\phi_{k-r}$ with
$\phi_{k-r} \in {\cal E}_{k-r}$.
The filtration of $\cal E$ implies $\phi_{k-r} \in {\cal E}_{k-r+1}$ so
$\nu^r\nu^{(k-r+1)}\phi_{k-r} \in {\cal U}_{k+1}$. The element
$\nu^{r+1} \nu^{(k-r)}\phi_{k-r}$ is also contained in ${\cal U}_{k+1}$.
By eqn. (\ref{nudef}) we get
$\nu^{(k-r+1)} = (1-(k-r)\nu)\nu^{(k-r)}$. Hence
\[
    \nu^r\nu^{(k-r)} \phi_{k-r} =
    \nu^r\nu^{(k-r+1)} \phi_{k-r} + (k-r) \nu^{r+1} \nu^{(k-r)} \phi_{k-r}
    \in {\cal U}_{k+1}
\]
which proves the inclusion ${\cal U}_k \subset {\cal U}_{k+1}$. It is
clear from the definition (\ref{dieUks}) that
$\nu {\cal U}_k \subset {\cal U}_{k+1} \subset {\cal E}[\nu]$
whence $\cal U$ is $\Bbb{C}[\nu]$-submodule of ${\cal E}[\nu]$.
Note that for fixed $y\in\CNull$ the function
$\pi(z) \mapsto \f{|\langle y,z \rangle|^4}{x^2}\in{\cal E}_2$ is
not contained in $\cal U$.
{\it ii.)} By Theorem (\ref{one}) the star product for arbitrary
$\Phi_k=\sum_{r=0}^k\nu^{k-r}\nu^{(r)}\phi_r \in {\cal U}_k$ and
$\Psi_l=\sum_{s=0}^l\nu^{l-s}\nu^{(s)}\psi_s \in {\cal U}_l$
with $\phi_r \in {\cal E}_r$ and $\psi_s \in {\cal E}_s$ is given by
\[
    \Phi_k \mustar \Psi_l =
    \sum_{r=0}^k\sum_{s=0}^l\sum_{t=0}^{min(r,s)}
    \frac{1}{t!} \nu^{k-r+l-s+t} \nu^{(r+s-t)}
    M_t^{(r,s)} (\phi_r, \psi_s),
\]
hence each summand is an element of ${\cal U}_{k+l}$.
This proves the second part.
\qed\\

We should now like to substitute the formal parameter $\nu$ for a fixed
nonzero real number $\alpha$. In the subalgebra $\cal U$ this
is well-defined because $\cal U$ only contains polynomials in $\nu$.
The kernel of the substitution homomorphism
\beq
    {\cal I}_\alpha:=\{\Phi(\nu)\in {\cal U}| \Phi(\alpha)=0\}
\eeq
will turn out to be a $\mustar$-ideal.
\begin{lemm}
    ${\cal I}_\alpha$ is a $\mustar$-ideal and the general form of an
    element $\Phi \in {\cal I}_\alpha \cap {\cal U}_k$ is
    \begin {enumerate}
    \item for $\alpha \not\in \{1, \frac{1}{2}, \frac{1}{3}, \ldots \}$
          \beq \label {alphaR}
            \Phi (\nu) = (\nu - \alpha) u_{k-1} (\nu)
            \quad \mbox { with } u_{k-1} (\nu) \in {\cal U}_{k-1}
          \eeq
    \item for $\alpha = \frac{1}{K}$ with $K \in \Bbb N\setminus \{0\}$
          for $k \le K$
          \beq
               \Phi (\nu) = (\nu - \frac{1}{K}) u_{k-1}(\nu)
               \label {alphakK}
          \eeq
          and for $k > K$
          \beq
               \Phi (\nu) = \nu^{(k)} \phi_k + \cdots +
               \nu^{k-K-1} \nu^{(k-K+1)} \phi_{k-K+1}
               + (\nu - \frac{1}{K}) u_{k-1} (\nu)
               \label {alphaK}
          \eeq
          with some $u_{k-1}(\nu) \in {\cal U}_{k-1}$ and
          $\phi_r \in {\cal E}_r$.
    \end {enumerate}
\end {lemm}

{\sc Proof :}
    It is clear that ${\cal I}_\alpha$ is a $\mustar$-ideal if the
    eqn. (\ref {alphaR}) resp. (\ref {alphakK}, \ref {alphaK})
    are valid because of the
    form of $\nu^{(k)}$ in eqn. (\ref {nudef}) and eqn. (\ref {StarEqn})
    for the $\mustar$-product. So we only have to prove eqn. (\ref {alphaR})
    and eqs. (\ref {alphakK}, \ref{alphaK}).

    Case {\it i.)}
    $\alpha \not\in \{ 1, \frac{1}{2}, \frac{1}{3},\ldots \}$
    We prove (\ref {alphaR}) by induction on k.
    With ${\cal U}_1 = \nu{\cal E}_0 + {\cal E}_1$ we get
    $\Phi (\nu) = c + \nu\sigma(A)$ for some $c \in \Bbb{C}$ and
    $\sigma (A) \in {\cal E}_1$. $\Phi (\alpha) = 0$ implies
    $A = - \frac{c}{\alpha} \mbox{\bf 1}$ where {\bf 1} is the unit matrix.
    So we have $\Phi(\nu) = c (1-\frac{\nu}{\alpha})$. This proves the
    case $k = 1$.

    Consider now the induction step $k\!-\!1 \to k$. Since
    ${\cal U}_k = \nu^{(k)} {\cal E}_k + \nu {\cal U}_{k-1}$ we write for
    $\Phi \in {\cal U}_k$
    \[
        \Phi (\nu) = \nu^{(k)} \phi_k + \nu u_{k-1} (\nu)
        \qquad \mbox {~ with~} \phi_k \in {\cal E}_k \mbox {~and~}
        u_{k-1} (\nu) \in {\cal U}_{k-1}.
    \]
    Then $\Phi (\alpha) = 0$ implies $\phi_k = -
    \frac{\alpha}{\alpha^{(k)}} u_{k-1} (\alpha) \in {\cal E}_{k-1}$
    because $\alpha^{(k)} \ne 0$.
    Writing $\nu^{(k)} = (1-(k-1)\nu)\nu^{(k-1)}$ we get
    \[
        \Phi (\nu) = \nu^{(k-1)} \phi_k + \nu \left(
        -(k-1)\nu^{(k-1)} \phi_k + u_{k-1} (\nu) \right)
        =: \nu^{(k-1)} \phi_k + \nu u_{k-1}'(\nu)
    \]
    with $u_{k-1}'(\nu) \in {\cal U}_{k-1}$.
    In $\Phi (\nu) = (\nu-\alpha) u_{k-1}'(\nu) + \nu^{(k-1)} \phi_k +
    \alpha u_{k-1}'(\nu)$ the first term vanishes at $\nu = \alpha$ so
    we have
    \[
        \left. \nu^{(k-1)} \phi_k + \alpha u_{k-1}'(\nu)
        \right|_{\nu=\alpha} = 0 .
    \]
    But $\nu^{(k-1)} \phi_k + \alpha u_{k-1}'(\nu)$ is an element in
    ${\cal U}_{k-1}$ so we can apply the induction and conclude
    \[
        \Phi (\nu) = (\nu - \alpha) u_{k-1}'(\nu) +
                     (\nu - \alpha) u_{k_2}'' (\nu)
    \]
    with some $u_{k-2}''(\nu) \in {\cal U}_{k-2}$. This proves the first
    part.

    Case {\it ii.)}
    Let $\alpha = \frac{1}{K}$ with $K \in \Bbb{N} \setminus \{0\}$
    and $\Phi (\nu) \in {\cal U}_k$. We have to consider
    the two cases $k \le K$ and $k > K$ separately:

    a.) For $k \le K$ we have $\left(\frac{1}{K}\right)^{(k)} \ne 0$ and
    we can apply the same arguments as in the first part, hence
    $\Phi (\nu) = (\nu - \frac{1}{K}) u_{k-1} (\nu)$ with
    $u_{k-1} (\nu) \in {\cal U}_{k-1}$.

    b.) For $k > K$ we have for $r \ge 1$
    $\left. \nu^{(K+r)}\right|_{\nu = \frac{1}{K}} = 0$ according to the
    definition (\ref{nudef}) of $\nu^{(k)}$.
    Hence in every $\Phi(\nu) \in {\cal U}_k$
    \[
        \Phi(\nu) = \nu^{(k)}\phi_k + \cdots + \nu^{k-K-1} \nu^{(K+1)}
                    \phi_{K+1} + \nu^{k-K} u_K (\nu)
    \]
    the first terms are automatically elements of
    ${\cal I}_{\frac{1}{K}}$ and $\Phi(\frac{1}{K}) = 0$ implies
    $u_K(\frac{1}{K}) = 0$. But this is an element of ${\cal U}_K$ and we can
    apply the case a.) Hence $u_K(\nu) = (\nu - \frac{1}{K}) u_{K-1}'(\nu)$
    which proves the second part. \qed \\

The quotient algebras ${\cal U}/{\cal I}_\alpha$ can now easily
be described:
\begin{theo}
   The quotient algebra ${\cal A}_\alpha:={\cal U}/{\cal I}_\alpha$
   is isomorphic to one of the following algebras:
   \begin{enumerate}
   \item For $\alpha$ not equal to one of the rational numbers
         $1,\f{1}{2},\f{1}{3},\ldots$ the algebra ${\cal A}_\alpha$
         is isomorphic to the vector space of representative
         functions $\cal E$ equipped with the multiplication $\ast_\alpha$
         defined by ($f\in{\cal E}_k,g\in{\cal E}_l;k,l\in\Bbb{N}$):
         \beq
             f\ast_\alpha g:=\sum_{r=0}^{min(k,l)}\f{\alpha^r}{r!}
             \f{\alpha^{(k+l-r)}}{\alpha^{(k)}\alpha^{(l)}}
             M^{(k,l)}_r(f,g)
         \eeq
         where the real number $\alpha^{(k)}$ is defined by the formula
         (\ref{nudef}) and $M_r^{(k,l)}(f,g)$ is given in (\ref{dieHs}).
   \item Let $\alpha$ be equal to $\f{1}{K}$ with $K$ a positive integer.
         Then ${\cal A}_\alpha$ is isomorphic to the finite-dimensional
         complex algebra of linear endomorphisms of $\Bbb{C}^N$ with
         $N:={n+K \choose K}$. The isomorphism
         is given by the map
         \beq \lb{finrep}
             A \mapsto  \f{\nu^{(K)}}{(1/K)^{(K)}}\sigma(A)
                                      \bmod {\cal I}_{1/K}
         \eeq
         where $\sigma(A)$ is the Berezin-Rawnsley symbol of the complex
         $N\times N$-matrix $A$ (see (\ref{BerRaw})). The matrix product
         $(AB)_{i_1\cdots i_K,j_1\cdots j_K}$ is given by
         $A_{i_1\cdots i_K,a_1\cdots a_K}$ $B_{a_1\cdots a_K,j_1\cdots j_K}$.
   \end{enumerate}
\end{theo}

{\sc Proof:} i) According to the preceding Lemma the ideal ${\cal I_\alpha}$
is equal to $(\nu-\alpha){\cal U}$ which amounts to simply substituting
$\nu=\alpha$ in (\ref{StarEqn}) which is obviously well-defined.

ii) According to the second part of the preceding Lemma we get
${\cal U}_k/{\cal I}_{1/K}={\cal U}_K/{\cal I}_{1/K}$ for $k\geq K$.
For ${\cal U}_K$ we may substitute $\nu$ for $1/K$ since $(1/K)^{(k)}\neq 0$
for $k\leq K$. Hence $\dim {\cal A}_{1/K}=\dim {\cal E}_K={n+K \choose K}^2$.
Moreover,
\beqa
  \lefteqn {\f{\nu^{(K)}}{(1/K)^{(K)}}\sigma(A)\mustar
          \f{\nu^{(K)}}{(1/K)^{(K)}}\sigma(B) \bmod {\cal I}_{1/K}} \\
                 & = &
  \sum_{r=0}^K\f{\nu^r}{r!}\f{\nu^{(2K-r)}}{(1/K)^{(K)}(1/K)^{(K)}}
               M^{(K,K)}_r(\sigma(A),\sigma(B)) \bmod {\cal I}_{1/K} \\
                 & = &
          \f{(1/K)^K}{K!}\f{\nu^{(K)}}{(1/K)^{(K)}(1/K)^{(K)}}
               M^{(K,K)}_K(\sigma(A),\sigma(B)) \bmod {\cal I}_{1/K} \\
                 & = &
        \f{\nu^{(K)}}{(1/K)^{(K)}}\sigma(AB) \bmod {\cal I}_{1/K} ,
\eeqa
since a simple calculation shows that
\[
  \pi^\ast(M^{(K,K)}_K(\sigma(A),\sigma(B)))(z)
    =
     K!K!\f{\tilde{\sigma}(AB)(z)}{x^K}.
\]
This proves the Theorem.
\qed\\

Remarks: Note that for each $A\in\goth{u}(n+1)$ the Berezin-Rawnsley
  symbol $\sigma(A)$ is contained in ${\cal E}_K$ and thus uniquely
  corresponds to a linear operator in $B({\cal H}^{(K)})$ which is mapped
  via the linear map (\ref{finrep}) to ${\cal A}_{1/K}$. By Lemma
  \ref{quantmom} it follows that this defines a
  {\em representation} of $\goth{u}(n+1)$ in $\Bbb{C}^N$, which is
  {\em irreducible}:
  in fact, by Lemma \ref{quantmom} we know that the momentum map
  $P(A)$ star-commutes with some function iff it Poisson commutes
  with that function which is only possible iff that function is constant
  since the unitary group acts transitively on $\CP^n$. Since the
  Poisson bracket with $P(A)$ obviously preserves each ${\cal E}_k$,
  it follows that this reasoning carries over to the quotient by
  ${\cal I}_{1/K}$.

  From the physical point of view the second part of the preceding Theorem
  can be viewed as follows:
         The classical phase space reduction of $\CNull$ to $\CP^n$
         was motivated by the $U(1)$-action on $\CNull$ which is
         just the flow of the classical isotropic $(n+1)$-dimensional
         harmonic
         oscillator with Hamiltonian $H = \frac{1}{2} x$ (The frequency
         and mass are normalised to 1).
         Passing from $\CNull$ to $\CP^n$ for a fixed value
         $\mu \in \Bbb{R}^-$ of the momentum mapping
         $J = -\frac{1}{2} x$ means classically that the
         harmonic oscillator is considered at a fixed energy $E = -\mu$.\\
         This is also true in the quantum mechanical case but now the
         energy is quantised. Only for the discrete values
         $\frac{1}{K}$, $K$ a positive integer, of the formal parameter
         $\nu = \frac{\lambda}{-2\mu}$ one obtains
         {\em finite dimensional} algebras of observables
         ${\cal U}/{\cal I}_{1/K}$ as one would physically expect for the
         {\em compact phase space} $\CP^n$ because the phase volume is
         finite and each state `occupies a phase cell of volume not smaller
         than $\hbar^n$' which results in a finite-dimensional
         Hilbert space for the quantum mechanical states.
         With $\lambda = 2\hbar$ the quantised energy is given by
         \[
            E_K = \hbar K
         \]
         (where the usual ground state energy $\frac{1}{2} \hbar (n+1)$
         is absent because of the Wick ordering). Note that in this
         interpretation the
         formal parameter $\lambda = 2\hbar$ is {\em not} quantised
         but the energy $E = -\mu$
         is. The dimensions of the operator algebras for a
         fixed $K\in \Bbb{N}$ correspond to the well-known degeneracy of
         the $K$th energy eigenvalue of the isotropic harmonic oszillator
         (see e.g. \cite[eqn XII.64]{Mes 61}).

\section{Other Examples}

In this section we briefly sketch how the programme mentioned in the
introduction applies to the deformation quantization of other phase spaces
which have been dealt with in the literature:

{\bf 1.} Consider complex $n+1$-space $\Bbb{C}^{n+1}$ as a symplectic
manifold in the usual manner, i.e. with symplectic form
$\omega= \f{{\bf i}}{2}\sum_{i=0}^n dz^i\wedge d\bar{z}^i$. The Wick product
(\ref{Wickusual}) then defines a star product on this space. It is
natural to consider the action of $\Bbb{C}^{n+1}$ on itself by translations.
Suppose that the smooth complex-valued function $F$ on $\Bbb{C}^{n+1}$ is
representative with respect to this group action, i.e. there is a finite
number $L$ of linearly independent smooth complex-valued
functions $F_1,\ldots,F_L$ on $\Bbb{C}^{n+1}$ such
that $F(z+v)=\sum_{a=1}^L\beta_a(v)F_a(z)$ with smooth complex valued
coefficients $\beta_a$. The same equation holds for each $F_a$ thus giving
rise to a coefficient matrix $\beta_{ab}(v)$. Since $\Bbb{C}^{n+1}$ is
abelian, $\beta_{ab}(v)$ commutes with each $\beta_{ab}(w)$ whence
all these matrices can simultaneously be transformed to Jordan normal form.
It can easily be seen that the generalized eigenvectors of $\beta_{ab}(v)$
are of the form $p(z)e^{b_iz^i+c_i\bar{z}^i}$ where $p$ is a complex-valued
polynomial function of $(z,\bar{z})$ and $b_i,c_j$ are complex numbers.
Conversely, it is easy to see that each linear combination of the functions
of this form
is indeed representative. Since this space of functions
is a subalgebra of the algebra of all smooth complex-valued functions
on $\Bbb{C}^{n+1}$ (with respect to pointwise multiplication) which
clearly is closed under complex conjugation, contains a unit element,
and separates points it is dense in $C^\infty(\Bbb{C}^{n+1})$ with
respect to the uniform topology on compacta thanks to the
Stone-Weierstrass Theorem. Note that the space of polynomials is
a dense proper $\Bbb{C}^{n+1}$ submodule of representative functions
which would be impossible for compact Lie groups (compare the proof of
Lemma \ref{Eistrep}). \\
Writing $e_{(a,b)}$ for the exponential function
$z\mapsto e^{a_iz^i+b_i\bar{z}^i}$ parametrized by two complex vectors
$a,b\in\Bbb{C}^{n+1}$ we easily obtain the following formula:
\beq
    e_{(a+\rho,b+\sigma)}\ast e_{(a'+\rho',b'+\sigma')}
           =
          e^{\lambda(a_i+\rho_i)(b'_i+\sigma'_i)}
           e_{(a+\rho,b+\sigma)}e_{(a'+\rho',b'+\sigma')}
\eeq
where $a,a',b,b',\rho,\rho',\sigma,\sigma'\in\Bbb{C}^{n+1}$. After
differentiating this formula a finite number of times with respect
to the components of $\rho,\rho',\sigma,\sigma'$ at
$\rho=\rho'=\sigma=\sigma'=0$ (which generates polynomial prefactors)
we see that the Wick-product of two representative functions
is again representative and entire analytic in the formal parmater
$\lambda$. Therefore the algebra $\cal U$ can be chosen to be the
subalgebra of $C^\infty(\Bbb{C}^{n+1})[[\lambda]]$ which consists of
polynomials in $\lambda$ with coefficients in the space of representative
functions. Substituting $\lambda$ for a real number $\alpha$ then is
straight-forward.

{\bf 2.} Consider now the $2n$-torus $T^{2n}:=S^1\times\cdots\times S^1$
($2n$ factors).
Let $(\vp^1,\ldots,\vp^{2n})=:\vp$, $0\leq\vp^i <1$, denote the standard
angle co-ordinates on $T^{2n}$ and take any nondegenerate
complex $2n\times 2n$-matrix $(\Lambda^{ij})$. The Moyal product may then
be defined by for two smooth, complex-valued functions $f,g$ on $T^{2n}$
as follows:
\beq \lb{Moyal}
    f\ast g:=\sum_{r=0}^\infty\f{(\lambda/2\pi{\bf i})^r}{r!}
                     \Lambda^{i_1j_1}\cdots\Lambda^{i_rj_r}
                     \f{\p^r f}{\p\vp^{i_1}\cdots\vp^{i_r}}
                     \f{\p^r g}{\p\vp^{j_1}\cdots\vp^{j_r}}
\eeq
It is known that this formula defines an associative deformation
for the pointwise multiplication in $C^\infty(T^{2n})$ (see e.g.
\cite{bayen}). The space of
representative functions for the torus action on itself is spanned by
the Fourier modes $T_k(\vp):=e^{2\pi{\bf i}k_i\vp^i}$, $k\in\Bbb{Z}^{2n}$:
indeed, it is clear
that the complex span of the Fourier modes consists of representative
functions. Since it is a subalgebra of $C^\infty(T^{2n})$ which is
closed under complex conjugation, contains a unit element, and separates
points it is a dense $T^{2n}$-submodule of the space of all representative
functions (by the Stone-Weierstrass Theorem) which has to be
equal to that space since $T^{2n}$ is compact (compare again the proof of
Lemma (\ref{Eistrep})). The Moyal product of two Fourier modes is then
simply be computed by
\beq \lb{MoyalFourier}
    T_k\ast T_{k'}=e^{2\pi{\bf i}\lambda\Lambda^{ij}k_ik'_j}T_kT_{k'}
\eeq
which again is an entire function in the formal parameter $\lambda$,
whence the choice of $\cal U$ and the substitution of the formal parameter
is completely analogous to the previous example. \\
Suppose now that the matrix $(\Lambda^{ij})$ is integral and the
greatest common divisor of the matrix elements is equal to $1$. Choosing
$\lambda=\f{1}{K}$, $K$ a positive integer, it is easily seen from the above
formula (\ref{MoyalFourier}) that the subspace ${\cal J}_{1/K}$ spanned
by all elements of the form $T_k-T_{k+Kk'}$ with $k,k'\in\Bbb{Z}^{2n}$
is a $\ast$-ideal in the $\ast$-algebra of representative functions.
It can be shown (see \cite{BHSS 91}) that the quotient algebra is a
simple complex algebra of dimension $K^{2n}$, which is related to the
geometric quantization on the $2n$-torus in the theta-bundle and its
tensor powers.

{\bf 3.} Let $\Bbb{D}:=\{v\in\Bbb{C}||v|^2<1\}$ be the Poincar\'e disk
in the complex plane. As we have indicated in the last section of
\cite{uns} we can use the projective representation of $\Bbb{D}$:
in $\Bbb{C}^2\setminus \{0\}$ consider the open subset defined by
the inequality $0<y:=|z^0|^2-|z^1|^2$ (the function $y$ was defined
with an erroneous sign in \cite{uns}).
The image of this open set under
the projection $\pi$ is an open set of $\CP^1$ which is holomorphically
diffeomorphic to $\Bbb{D}$ via $\pi(z)\mapsto v:=\f{z^1}{z^0}$. In their
article \cite{CGR III} Cahen, Gutt, and Rawnsley have considered
the following functions $f_{p,q}(v)$ ($p,q\in\Bbb{N}$) on $\Bbb{D}$
\beq
  f_{p,q}(v):=v^p(\f{\bar{v}}{1-|v|^2})^q.
\eeq
Their pull-back to $\Bbb{C}^2\setminus \{0\}$ is simply given by
$(\pi^*f_{p,q})(z)=\f{(z^1/z^0)^p(z^0\bar{z}^1)^q}{y^q}$.
In \cite{uns} we have defined a star product
on $\Bbb{D}$ by essentially replacing $x$ by $y$ and the operators
$M_r$ by
\beq
    \check{M}_r(G,H):=y^rg^{i_1j_1}\cdots g^{i_rj_r}
                      \f{\p^rG}{\p z^{i_1}\cdots\p z^{i_r}}
                      \f{\p^rG}{\p \bar{z}^{j_1}\cdots\p \bar{z}^{j_r}}
\eeq
with $g^{ij}:={\rm diag}(1,-1)$.
Observing that eqs (\ref{oneeins}) and (\ref{StarEqn}) remain valid
for arbitrary smooth complex-valued functions when the upper bound
of the sum is $\infty$ and that these formulas pass to the noncompact
case with the above replacements and adapting the sign and ordering
conventions of \cite{CGR III} we obtain the following star product:
\beqa
    { (f_{p,q}\ast f_{r,s})(\pi(z)) }
        & = &
                 \sum_{m=0}^\infty\f{(-\nu)^m}{m!}
                         \f{(-\nu)^{(q+s-m)}}{(-\nu)^{(q)}(-\nu)^{(s)}}
                            y^{m-q-s}
                         g^{i_1j_1}\cdots g^{i_mj_m} \\
        &   & \;\; \times \;\;
              \f{\p^r(y^q\pi^*f_{p,q})}
              {\p \bar{z}^{i_1}\cdots\p \bar{z}^{i_m}}(z)
              \f{\p^r(y^s\pi^*f_{r,s})}
              {\p z^{j_1}\cdots\p z^{j_m}}(z) \\
        & = & \sum_{m=0}^\infty\f{(-\nu)^m}{m!}
              \f{(-\nu)^{(q+s-m)}}{(-\nu)^{(q)}(-\nu)^{(s)}}
              y^{m-q-s}(-1)^m \\
        &   & \;\; \times \;\;
              \f{\p^m(((z^1/z^0)^p(z^0\bar{z^1})^q)}{\p (\bar{z^1})^m}(z)
              \f{\p^m(((z^1/z^0)^r(z^0\bar{z^1})^s)}{\p (z^1)^m}(z)   \\
        & = & \sum_{m=0}^{\min (q,r)}\f{\nu^m}{m!}
              \f{(-\nu)^{(q+s-m)}}{(-\nu)^{(q)}(-\nu)^{(s)}}
              \f{q!}{(q-m)!}\f{r!}{(r-m)!}
              f_{p+r-m,q+s-m}(\pi(z))
\eeqa
which reproduces the result of \cite{CGR III}. In \cite{CGR IV} star products
were computed on more general bounded symmetric domains.

\vspace{0.5cm}

\noindent
{\large\bf Acknowledgment}

\vspace{2mm}

\noindent The authors would like to thank D.~Arnal, J.~Hoppe, M.~Masmoudi
and W.~Soergel for valuable discussions.

\end{document}